\documentclass[12pt, reqno, dvipdfmx]{amsart}
\usepackage{preemble}

\begin{document}

\title[Non-noise sensitivity for word hyperbolic groups]
      {Non-noise sensitivity for word hyperbolic groups}
\author{Ryokichi Tanaka}
\address{Department of Mathematics, 
      Kyoto University, Kyoto 606-8502 JAPAN}
\email{rtanaka@math.kyoto-u.ac.jp}
\date{\today}

\maketitle

\begin{abstract}
We show that non-elementary random walks on word hyperbolic groups with finite first moment are not noise sensitive in a strong sense for small noise parameters.
\end{abstract}

\section{Introduction}

Let $\G$ be a countable group, and $\m$ be a probability measure on it.
The main interest is in the case when $\G$ is finitely generated with a finite set of generators $S$ and $\m$ is the uniform distribution on $S$,
but we will also consider the case when $\m$ has an unbounded support.
The $\m$-random walk 
starting from the identity $o$
is defined by a product of independent sequence of increments with the identical distribution $\m$.
The noise sensitivity question concerning a $\m$-random walk on $\G$ asks the following:
If we choose some real parameter $\rho \in (0, 1)$ and replace each increment by an independent sample with the same law $\m$ with probability $\rho$ or retain it with probability $1-\rho$, independently,
then is the resulting random walk asymptotically independent of the original one?

More precisely, the $\m$-random walk $\{w_n\}_{n=0}^\infty$ starting from $o$ is defined by $w_n=x_1\cdots x_n$ for an independent sequence $\{x_n\}_{n=1}^\infty$ with the identical distribution $\m$ and $w_0:=o$.
Let $\m_n$ denote the distribution of $w_n$, which is the $n$-fold convolutions of $\m$.
For each $\rho \in [0, 1]$,
let us consider the $\pi^\rho$-random walk $\{(w_n^{(1)}, w_n^{(2)})\}_{n=0}^\infty$ on $\G \times \G$ defined by
\[
\pi^\rho:=\rho \m\times \m+(1-\rho)\m_{\diag} \quad \text{on $\G \times \G$},
\]
where $\m \times \m$ denotes the product measure and 
$\m_{\diag}((x, y))=\m(x)$ if $x=y$ and $0$ otherwise.
For any two probability measures $\n_1$ and $\n_2$ on a countable set $X$,
the total-variation distance is defined by
\[
\|\n_1-\n_2\|_{\TV}:=\sup_{A \subset X}|\n_1(A)-\n_2(A)|.
\]

\begin{definition}
The $\m$-random walk on $\G$ is called $\ell^1$-\textit{noise sensitive} if for all $0<\rho < 1$,
\[
\|\pi^\rho_n-\m_n \times \m_n\|_{\TV} \to 0 \quad \text{as $n \to \infty$}.
\]
\end{definition}

The $\ell^1$-noise sensitivity was introduced by Benjamini and Brieussel \cite[Definition 2.1]{BenjaminiBrieussel}.
There it has been shown that this notion in general highly depends not only on $\G$ but also on $\m$.
Among others, they have proved that if $\G$ admits a surjective homomorphism onto the infinite cyclic group $\Z$ and the support of $\m$ is a finite set of generators, then a $\m$-random walk on $\G$ is \textit{not} $\ell^1$-nose sensitive.
Moreover, if $(\G, \m)$ is \textit{non-Liouville}, i.e., there exists a non-constant bounded $\m$-harmonic function on $\G$,
then a $\m$-random walk on $\G$ is \textit{not} $\ell^1$-noise sensitive \cite[Theorem 1.1]{BenjaminiBrieussel}.
It is believed that these two properties are the only possible obstructions for the $\ell^1$-noise sensitivity for random walks on groups.
We show that non-elementary word hyperbolic groups with large class of $\m$ reveal a strong negation of the $\ell^1$-noise sensitivity if $\rho$ is small enough.
This also offers in the non-Liouville setting a way to show that random walks are not $\ell^1$-noise sensitive in some refined sense for a class of groups possibly without non-trivial homomorphisms onto $\Z$.

Let $\G$ be a word hyperbolic group.
A probability measure $\m$ on $\G$ is called \textit{non-elementary} if the support generates a non-elementary subgroup $\gr(\m)$ as a group.
In this setting, it is equivalent to say that $\gr(\m)$ contains a free group of rank greater than one (and $\G$ is necessarily non-elementary).
Furthermore we say that $\m$ has \textit{finite first moment} if 
\[
\sum_{x \in \G}|x|\m(x)<\infty,
\]
for some (equivalently, every) word norm $|\cdot|$.
First we note that all $\m$-random walk on $\G$ for a non-elementary $\m$ is not $\ell^1$-noise sensitive without any moment condition.

\begin{theorem}\label{Thm:1}
Let $\G$ be a word hyperbolic group and $\m$ be a non-elementary probability measure on $\G$.
For all $0\le \rho < 1$, we have
\[
\liminf_{n \to \infty}\|\pi^\rho_n-\m_n \times \m_n\|_{\TV}>0.
\]
\end{theorem}

This in fact follows from the proof of \cite[Theorem 4.1]{BenjaminiBrieussel} due to the non-Liouville property.
We provide a proof in this setting to illustrate our approach.
Second we show that under the finite first moment condition if $\rho$ is close enough to $0$, then the distribution of a $\pi^\rho$-random walk and the joint distribution of independent copies of two $\m$-random walks are mutually singular at infinity.

\begin{theorem}\label{Thm:2}
Let $\G$ be a word hyperbolic group and $\m$ be a non-elementary probability measure with finite first moment on $\G$.
There exists some $0<\rho_\ast\le 1$ such that for all $0<\rho<\rho_\ast$, we have
\[
\|\pi^\rho_n-\m_n\times \m_n\|_\TV \to 1 \quad \text{as $n \to \infty$}.
\]
\end{theorem}

In the case when $\m$ is not a non-elementary probability measure,
these statements are no longer true.
Indeed, it suffices to consider elementary word hyperbolic groups, which are either finite groups or contain $\Z$ as a finite index subgroup (see e.g., \cite{Gromovhyperbolic} and \cite{GhysdelaHarpe} for background).
If $\G$ is a finite group, then for every probability measure $\m$ on it, a $\m$-random walk on $\G$ is $\ell^1$-noise sensitive. Indeed, $\pi^\rho$ and $\m\times \m$ have the same support on $\G \times \G$ and the distributions of corresponding random walks with the same initial state tend to a common stationary distribution for time with the same parity (cf.\ \cite[Proposition 5.1]{BenjaminiBrieussel}).
Benjamini and Brieussel have shown that on the infinite dihedral group, a lazy simple random walk on a Cayley graph is $\ell^1$-noise sensitive \cite[Theorem 1.4]{BenjaminiBrieussel}.

As it is mentioned above, if $\G$ has a surjective homomorphism onto $\Z$,
then a $\m$-random walk can not be $\ell^1$-noise sensitive.
This follows by computing covariance of the random walk on each factor of $\Z^2$ as the image of product of homomorphisms and the classical central limit theorem.
The method thus works for $\m$ with finite second moment.
We note, however, that so far this and the non-Liouville property have been the only known ways to disprove $\ell^1$-noise sensitivity.
See \cite{BenjaminiBrieussel} and \cite[Section 3.3.4]{KalaiICM2018} for discussions concerning the subject of matters, various interesting notions of noise sensitivity for random walks on groups, questions and conjectures.

The proofs of Theorems \ref{Thm:1} and \ref{Thm:2} rely on the boundary of word hyperbolic groups,
in particular, the fact that a $\m$-boundary or the Poisson boundary for $(\G, \m)$ is realized on a topological boundary of the group.
In this setting, we show that if $h(\pi^\rho)\neq h(\pi^{\rho'})$ for $0 \le \rho, \rho' \le 1$, 
then
\begin{equation}\label{Eq:intro_pi^rho}
\|\pi^\rho_n-\pi^{\rho'}_n\|_\TV \to 1 \quad \text{as $n \to \infty$},
\end{equation}
where $h(\pi^\rho)$ is the asymptotic entropy for a $\pi^\rho$-random walk
(see Section \ref{Sec:preliminary} for the definition, and Remark \ref{Rem:entropy}).
It is proven by showing that the harmonic measure $\n_{\pi^\rho}$ on the product of boundaries $(\partial \G)^2$ is \textit{exact dimensional} with a natural (quasi-) metric, i.e.,
\[
\frac{\log \n_{\pi^\rho}(B(\xib, r))}{\log r} \to \frac{h(\pi^\rho)}{l} \quad \text{as $r \to 0$ for $\n_{\pi^\rho}$-almost every $\xib$ in $(\partial \G)^2$},
\]
where $l$ stands for the drift defined by a product metric in $\G \times \G$
(see Theorem \ref{Thm:dim} for the precise statement).
The proof follows the methods in \cite{Tdim}, adapting to a product of word hyperbolic groups.
By \eqref{Eq:intro_pi^rho}, together with the continuity of $h(\pi^\rho)$ in $\rho \in [0, 1]$ (Corollary \ref{Cor:conti}),
as it is established by a general result of Erschler and Kaimanovich \cite{ErschlerKaimanovich},
we show that there exists some $0<\rho_\ast\le 1$ with $h(\pi^\rho)<h(\pi^1)$ for all $0<\rho<\rho_\ast$, 
deducing Theorem \ref{Thm:2}.

It would be interesting to determine $\rho_\ast$ in Theorem \ref{Thm:2}.
For example, 
in the case when $\m$ is a uniform distribution on a finite set of size $m \ge 2$, freely generating a free \textit{semi}-group,
the asymptotic entropy is explicitly computed as 
\[
h(\pi^\rho)=\log m-\(1-\frac{m-1}{m}\rho\)\log\(1-\frac{m-1}{m}\rho\)-\frac{m-1}{m}\rho\log \frac{\rho}{m} \quad \text{for $0 \le \rho \le 1$},
\] 
and $h(\pi^\rho)<h(\pi^1)=2\log m$ if $\rho<1$.
This implies that $\rho_\ast=1$ in the special case.
It might be expected that $\rho_\ast=1$ in many cases, however, we do not know how to show this in general.

Organization of this paper is the following:
In Section \ref{Sec:preliminary} we review known facts and tools on random walks and word hyperbolic groups,
in Section \ref{Sec:dim} we show that the harmonic measure for $\pi$ whose marginals are a common non-elementary probability measure $\m$ on $\G$ with finite first moment is exact dimensional in Theorem \ref{Thm:dim}, in Section \ref{Sec:conti} the continuity of asymptotic entropy in the parameter is established in Corollary \ref{Cor:conti}, following \cite{ErschlerKaimanovich},
in Section \ref{Sec:proofs} we show Theorems \ref{Thm:1} and \ref{Thm:2},
and in Appendix \ref{Appendix}, we give a review concerning Poisson boundary for random walks and a proof on continuity of asymptotic entropy for the sake of convenience, mainly for an expository purpose.

\subsection*{Notations}
For a real valued function $f$ on the set of non-negative integers, we write $f(n)=o(n)$ if $|f(n)|/n \to 0$ as $n \to \infty$ and $f(n)=O(n)$ if there exists a constant $C$ such that $|f(n)|\le C n$ for all large enough $n$.
For a set $A$, we denote by $\#A$ the cardinality, and by $A^{\sf c}$ its complement set.
The set of non-negative integers is written as $\Z_+=\{0, 1, 2, \dots\}$.
We define $0\log 0:=0$.

\section{Preliminary}\label{Sec:preliminary}

\subsection{Random walks on word hyperbolic groups and their products}

For a countable group $\G$ and a probability measure $\m$ on it,
the \textit{asymptotic entropy} of a $\m$-random walk is defined by
\[
h(\m)=\lim_{n \to \infty}\frac{1}{n}H(\m_n),
\]
where $H(\m):=-\sum_{x \in \G}\m(x)\log \m(x)$, the \textit{Shannon entropy} for probability measures $\m$ on $\G$,
the limit exists by sub-additivity of $n\mapsto H(\m_n)$ and is finite if $H(\m)<\infty$.

Let $\G$ be a word hyperbolic group.
For a probability measure $\m$ on $\G$, let $\supp \m$ denote the support and we assume that the group $\gr(\m)$ generated by $\supp \m$ as a group is a non-elementary subgroup.
We fix a left-invariant word metric $d$ associated with some finite set of generators $S$ closed under inversion $s \mapsto s^{-1}$.
The following discussion does not depend on the choice of $S$.
We denote the associated distance function from the identity $o$ by $|x|:=d(o, x)$.
If $\m$ has finite first moment, then $H(\m)<\infty$ (cf.\ \cite[Section VII, B]{Derriennic}).

Let us consider any probability measure $\pi$ on $\G \times \G$ such that the push-forward of $\pi$ on each factor is a fixed $\m$ on $\G$. 
Let $(\O, \Fc, \Pb)$ be the probability measure space, where $\O=(\G \times \G)^{\Z_+}$, $\Fc$ is the $\s$-field generated by the cylinder sets and $\Pb$ is the distribution of the $\pi$-random walk $\{\wb_n\}_{n=0}^\infty$ starting from the identity on $\G \times \G$.
The expectation relative to $\Pb$ is denoted by $\Eb$.
Let $\pi_n$ be the distribution of $\wb_n$.
Note that for $\wb_n=(w_n^{(1)}, w_n^{(2)})$, each $\{w_n^{(i)}\}_{n=0}^\infty$ gives the $\m$-random walk on $\G$ starting from $o$ for $i=1, 2$. 
In this setting, we have for all $n \ge 0$,
\[
H(\m_n)\le H(\pi_n)\le 2 H(\m_n),
\]
and the asymptotic entropy $h(\pi)$ of a $\pi$-random walk is finite since $H(\m)<\infty$, and $h(\pi)$ is positive
since $h(\pi) \ge h(\m)>0$ and $\gr(\m)$ is a non-elementary subgroup in $\G$.
We define the metric 
\[
d_\times((x_1, x_2), (y_1, y_2)):=\max\{d(x_i, y_i), i=1, 2\} \quad \text{for $(x_1, x_2), (y_1, y_2) \in \G\times \G$}.
\]
Suppose that $\m$ has finite first moment.
Then $\pi$ has finite first moment relative to the distance function $d_\times(o, \cdot\,)$.
The \textit{drift} is defined by the limit
\[
l:=\lim_{n \to \infty}\frac{1}{n}\Eb d_\times(o, \wb_n),
\]
where the limit exists by sub-additivity $n \mapsto \Eb d_\times(o, \wb_n)$ and is finite.
The value $l$ coincides with the drift of a $\m$-random walk relative to $|\cdot|$, i.e.,
\begin{equation}\label{Eq:drift}
l=\lim_{n \to \infty}\frac{1}{n} |w_n^{(i)}| \quad \text{for $i=1, 2$ and for $\Pb$-almost every $\omega$ in $\Omega$},
\end{equation}
and also in $L^1(\Omega, \Fc, \Pb)$ by the Kingman subadditive ergodic theorem, 
and $l>0$ since $\gr(\m)$ is non-elementary \cite[Section 7.3]{Kaimanovich-hyperbolic}.

\subsection{Word hyperbolic groups}
We refer to \cite{Gromovhyperbolic} and \cite{GhysdelaHarpe} for background.
Let $\partial \G$ denote the (Gromov) boundary and we endow $\G \cup \partial \G$ with the natural topology which is compact and metrizable.
Letting $(x|y)_o$ be the Gromov product for $x, y \in \G \cup \partial \G$ based at $o$,
we define the quasi-metric in $\partial \G$ by
\[
q(\x, \y):=e^{-(\x|\y)_o} \quad \text{for $\x, \y \in \partial \G$}.
\]
Note that $q$ satisfies that $q(\x, \y)=0$ if and only if $\x=\y$, $q(\x, \y)=q(\y, \x)$ for $\x, \y \in \partial \G$ and the triangle inequality holds up to a multiplicative constant independent of the points.
It is known that $q$ is bi-H\"older equivalent to a genuine metric in $\partial \G$ yielding the original topology.
We work with the quasi-metric $q$ to define balls: Let 
\[
B(\x, r):=\big\{\y \in \partial \G \ : \ q(\x, \y)\le r\big\} \quad \text{for $\x \in \partial \G$ and for real $r \ge 0$}.
\]
For any positive real $R>0$ and $x \in \G$,
the \textit{shadow} is defined by
\[
\Oc(x, R):=\big\{\x \in \partial \G \ : \ (x|\x)_o \ge |x|-R\big\}.
\]
By the hyperbolicity of geodesic metric space $(\G, d)$,
for each fixed $T>0$,
there exist constants $R_0$, $C>0$ such that for all $R >R_0$, all $\x \in \partial \G$ and all $x \in \G$ in a $T$-neighborhood of a geodesic ray from $o$ to $\x$, 
we have
\begin{equation}\label{Eq:shadow_ball}
B(\x, C^{-1}e^{-|x|+R})\subset \Oc(x, R) \subset B(\x, C e^{-|x|+R}).
\end{equation}
In the product space $(\partial \G)^2$, 
we define
\[
q_\times((\x_1, \x_2), (\y_1, \y_2)):=\max\{q(\x_i, \y_i), i=1, 2\} \quad \text{for $(\x_1, \x_2), (\y_1, \y_2) \in (\partial \G)^2$}.
\]
By the definition, the ball of radius $r$ centered at $\xib$ in $(\partial \G)^2$ relative to $q_\times$ is obtained by
\begin{equation}\label{Eq:prod_ball}
B(\xib, r)=B(\x_1, r) \times B(\x_2, r) \quad \text{where $\xib=(\x_1, \x_2)$}.
\end{equation}

\section{The dimension of harmonic measure}\label{Sec:dim}

The $\m$-random walk $\{w_n\}_{n=0}^\infty$ on a word hyperbolic group $\G$ converges to a point $w_\infty$ in $\partial \G$ almost surely as $n \to \infty$ in $\G \cup \partial \G$ if $\gr(\m)$ is non-elementary \cite[Theorem 7.6]{Kaimanovich-hyperbolic}.
This implies that
the $\pi$-random walk $\{\wb_n\}_{n=0}^\infty$ on $\G \times \G$ converges to a point $\wb_\infty:=(w_\infty^{(1)}, w_\infty^{(2)})$ in $(\partial \G)^2$ almost surely as $n \to \infty$ in the product space $(\G \cup \partial \G)^2$, where $\wb_n=(w_n^{(1)}, w_n^{(2)})$ and each $w_n^{(i)}$ tends to $w_\infty^{(i)}$ almost surely for $i=1, 2$.
Let $\n_\m$ and $\n_\pi$ denote the limiting distribution of $w_\infty$ on $\partial \G$ and that of $(w_\infty^{(1)}, w_\infty^{(2)})$ on $(\partial \G)^2$, respectively.
We call $\n_\m$ and $\n_\pi$ \textit{harmonic measures} on $\partial \G$ and on $(\partial \G)^2$, respectively.
Note that the push-forward of $\n_\pi$ on each factor $\partial \G$ coincides with $\n_\m$.
The harmonic measure $\n_\pi$ (resp.\ $\n_\m$) is \textit{$\pi$-stationary} (resp.\ $\m$-stationary), i.e.,
\begin{equation}\label{Eq:pi-stationary}
\pi\ast \n_\pi=\n_\pi \quad \text{on $(\partial \G)^2$},
\end{equation}
where $\pi\ast\n_\pi:=\sum_{\xb \in \G \times \G}\pi(\xb)\xb\n_{\pi}$ and $\xb\n_\pi=\n_\pi\circ \xb^{-1}$,
and similarly,
\begin{equation}\label{Eq:mu-stationary}
\m\ast \n_\m=\n_\m \quad \text{on $\partial \G$}.
\end{equation}
The $\n_\pi$ and $\n_\m$ are unique such measures satisfying \eqref{Eq:pi-stationary} and \eqref{Eq:mu-stationary}, respectively \cite[cf.\ Theorem 2.4]{Kaimanovich-hyperbolic}.
Concerning more on background, see \cite{Kaimanovich-hyperbolic}.

For the harmonic measure $\n_\pi$ for the $\pi$-random walk, we show the following:

\begin{theorem}\label{Thm:dim}
Let $\G$ be a word hyperbolic group and $\m$ be a non-elementary probability measure on $\G$ with finite first moment.
If $\pi$ is a probability measure on $\G\times \G$ such that the push-forward of $\pi$ on each factor $\G$ is $\m$,
then
the corresponding harmonic measure $\n_\pi$ on $(\partial \G)^2$ is exact dimensional, i.e.,
\[
\lim_{r\to 0}\frac{\log \n_\pi(B(\xib,r))}{\log r}=\frac{h(\pi)}{l} \quad \text{for $\n_\pi$-almost every $\xib$ in $(\partial \G)^2$},
\]
where $h(\pi)$ is the asymptotic entropy, $l$ is the drift relative to $d_\times$ for the $\pi$-random walk on $\G\times \G$ and the ball $B(\xib, r)$ is defined by $q_\times$ in $(\partial \G)^2$.
\end{theorem}
\def\H{{\rm H}}
Let us define the (upper) Hausdorff dimension of $\n_\pi$ by
\[
\dim_\H \n_\pi=\inf\big\{ \dim_\H E \ : \ \text{$\n_\pi(E)=1$ and $E$ is Borel}\ \big\},
\]
where $\dim_\H E$ stands for the Hausdorff dimension of $E$ relative to $q_\times$ in $(\partial \G)^2$.
Theorem \ref{Thm:dim} together with the Frostman-type lemma (cf.\ \cite[Section 2.2]{Tdim}) shows the following:

\begin{corollary}
In the setting of Theorem \ref{Thm:dim},
we have
\[
\dim_\H \n_\pi=\frac{h(\pi)}{l}.
\]
\end{corollary}

Let us keep the same setting and notations as in Theorem \ref{Thm:dim} throughout this section.
We use the following \textit{ray approximation} of a $\m$-random walk on a word hyperbolic group $\G$:
Let $\Pi$ be a Borel measurable map from $\partial \G$ to the space $\Pc$ of unit speed geodesic rays from $o$ in $(\G, d)$ endowed with the topology of convergence on compact sets.
(Here $\Pi$ is defined as a Borel measurable map by choosing a total order on a fixed set of generators and the lexicographical minimal geodesic for each point in the boundary.)
Letting
\[
\g_\x=\Pi(\x) \quad \text{for $\x \in \partial \G$},
\]
we have
\begin{equation}\label{Eq:ray}
d(w_n(\o), \g_{w_\infty(\o)}(l n))=o(n) \quad \text{for $\Pb$-almost every $\o$ in $\O$},
\end{equation}
\cite[Section 7.4]{Kaimanovich-hyperbolic}, where one should write instead $\lfloor ln \rfloor$ (the integer part of $ln$) here and below, however, we keep ``$ln$" for the sake of simplicity.
Let us define the map
\[
\Pi_\times: (\partial \G)^2 \to \Pc \times \Pc, \quad (\x, \y)\mapsto (\g_\x, \g_\y),
\]
as a Borel measurable map.
For a $\pi$-random walk $\{\wb_n\}_{n=0}^\infty$ on $\G \times \G$, we have by \eqref{Eq:ray},
\begin{equation}\label{Eq:ray2}
d_\times(\wb_n(\o), \Pi_\times(\wb_\infty(\o))(ln))=o(n) \quad \text{for $\Pb$-almost every $\o$ in $\O$}.
\end{equation}
Recall that the Shannon theorem for random walks:
\begin{equation}\label{Eq:Shannon}
h(\pi)=\lim_{n \to \infty}-\frac{1}{n}\log \pi_n(\wb_n(\o)) \quad \text{for $\Pb$-almost every $\o$ in $\O$},
\end{equation}
which follows from the Kingman subadditive ergodic theorem (\cite[Theorem 2.1]{KaimanovichVershik} and \cite[Section IV]{Derriennic80} where Y.\ Derriennic attributes to an observation by J.\ P.\ Conze).

First we show the dimension upper bound in the claim.

\begin{lemma}\label{Lem:upperbound}
For $\n_\pi$-almost all $\xib$ in $(\partial\G)^2$,
\[
\limsup_{r \to 0}\frac{\log \n_\pi(B(\xib, r))}{\log r} \le \frac{h(\pi)}{l}.
\]
\end{lemma}

\proof
For all $\e>0$ and all interval $I$ in $[0, \infty)\cap \Z$,
let
\begin{equation*}
A_{\e, I}:=
\bigcap_{n \in I}\Bigg\{ \o \in \O \ : \ 
\begin{aligned}
& |d(o, w_n^{(i)}(\o))-l n| \le \e n, \ |w_{n}^{(i)}(\o)^{-1}w_{n+1}^{(i)}(\o)|\le \e n \ \text{for $i=1, 2$}\\
& \text{and} \ \pi_n(\wb_n(\o)) \ge e^{-n(h(\pi)+\e)}
\end{aligned}
\Bigg\}.
\end{equation*}
Since $\mu$ has finite first moment, $|w_{n}^{(i)}(\o)^{-1}w_{n+1}^{(i)}(\o)|\le \e n$ for $i=1, 2$ for all large $n$ for $\Pb$-almost every $\o$ in $\Omega$, and
by \eqref{Eq:drift} and \eqref{Eq:Shannon},
there exists an $N_\e$ such that $\Pb(A_{\e, [N_\e, \infty)}) \ge 1-\e$.
Let $A:=A_{\e, [N_\e, \infty)}$.
For each $\o \in \O$,
let
\[
C_n(\o):=\big\{\y \in \O \ : \ \wb_n(\y)=\wb_n(\o)\big\} \quad \text{for $n \ge 0$},
\]
which defines the event where a $\pi$-random walk after time $n$ is $\wb_n(\o)$.
The conditional probabilities satisfy that
\begin{equation}\label{Eq:PA}
\liminf_{n \to \infty}\Pb(A\mid C_n(\o))>0 \quad \text{for $\Pb$-almost every $\o \in A$}.
\end{equation}
Indeed, 
letting $A_{[N, n)}:=A_{\e, [N_\e, n)}$ and $A_{[n, \infty)}:=A_{\e, [n, \infty)}$ for simplicity of notation,
we have $A=A_{[N, n)}\cap A_{[n, \infty)}$ for $n> N$,
and by the Markov property of the $\pi$-random walk,
\[
\Pb(A\mid C_n(\o))=\Pb(A_{[N, n)}\mid C_n(\o))\Pb(A_{[n, \infty)}\mid C_n(\o)).
\]
Let $\sigma(w_n, w_{n+1}, \dots)$ denote the $\sigma$-algebra generated by $w_n, w_{n+1}, \dots$.
Since for $\Pb$-almost every $\o \in A=A_{[N, n)}\cap A_{[n, \infty)}$,
\begin{align*}
\Pb(A_{[N, n)}\mid C_n(\o))
&=\Pb(A_{[N, n)}\mid \s(w_n, w_{n+1}, \dots))(\o)\\
&=\Pb(A\mid\s(w_n, w_{n+1}, \dots))(\o),
\end{align*}
we have by the bounded martingale convergence theorem,
\[
\Pb(A_{[N, n)}\mid C_n(\o)) \to \Pb(A\mid \Tc)(\o) \quad \text{for $\Pb$-almost every $\o \in A$},
\]
as $n \to \infty$, where $\Tc:=\bigcap_{n=0}^\infty \s(w_n, w_{n+1}, \dots)$.
Furthermore, since for $\Pb$-almost every $\o \in A$,
\begin{align*}
\Pb(A_{[n, \infty)}\mid C_n(\o))
&=\Pb(A_{[n, \infty)}\mid \s(w_0, w_1, \dots, w_n))(\o)\\
&=\Pb(A\mid\s(w_0, w_1, \dots, w_n))(\o),
\end{align*}
we have
\[
\Pb(A_{[n, \infty)}\mid C_n(\o))\to \1_{A}(\o) \quad \text{for $\Pb$-almost every $\o \in A$},
\]
as $n \to \infty$.
Note that $\Pb(A\mid \Tc)(\o)>0$ for $\Pb$-almost every $\o \in A$.
Indeed, letting $A_{>0}:=\{\omega \in \Omega \ : \ \Pb(A\mid \Tc)(\omega)>0\}$,
we have $\Pb(A\mid \Tc)=\Pb(A\mid \Tc)\1_{A_{>0}}$ almost everywhere in $\Pb$, whence integrating both sides yields $\Pb(A)=\Pb(A \cap A_{>0})$.
Thus we obtain \eqref{Eq:PA}.

For $\Pb$-almost every $\o \in A$,
we have for each $i=1, 2$,
\[
|d(o, w_n^{(i)}(\o))- ln|\le \e n \quad \text{and} \quad d(w_{n}^{(i)}(o), w_{n+1}^{(i)}(\o))\le \e n \quad \text{for all $n \ge N_\e$},
\]
whence $w_\infty^{(i)}$ is defined and  
\[
(w_n^{(i)}(\o)\mid w_\infty^{(i)}(\o))_o \ge (l-2 \e) n- R \quad \text{for all $n \ge N_\e$},
\]
for a constant $R \ge 0$ independent of $\o$ or $n$ (cf.\ \cite[Section 7.2]{Kaimanovich-hyperbolic}).
For $\Pb$-almost every $\y \in A\cap C_n(\o)$,
since $\wb_n(\y)=\wb_n(\o)$, 
by the $\delta$-hyperbolicity we have 
\[
(w_\infty^{(i)}(\eta)\mid w_\infty^{(i)}(\o))_o \ge (l-2 \e) n- R-\delta \quad \text{for each $i=1, 2$},
\]
and thus we obtain by \eqref{Eq:prod_ball}, for $\Pb$-almost every $\o \in A$,
\[
\wb_\infty(\y) \in B\(\wb_\infty(\o), C e^{-(l-2\e)n}\) \quad \text{for $\Pb$-almost every $\y \in A\cap C_n(\o)$},
\]
where $C=e^{R+\delta}$ is a positive constant depending only on the metric of the group $\G$.
Therefore for $\Pb$-almost every $\o \in A$,
\[
\Pb(A\cap C_n(\o)) \le \n_\pi\(B\(\wb_\infty(\o), C e^{-(l-2\e)n}\)\).
\]
Moreover, by the definition of $A$, we have $\Pb(C_n(\o))=\pi_n(\wb_n(\o)) \ge e^{-n(h(\pi)+\e)}$ for all $n \ge N_\e$.
Invoking \eqref{Eq:PA}, we obtain
\[
\limsup_{n \to \infty}\frac{\log \n_\pi\(B\(\wb_\infty(\o), C e^{-(l-2\e)n}\)\)}{-n} \le h(\pi)+\e
\quad \text{for $\Pb$-almost every $\o \in A$}.
\]
Noting that $r_n:=C e^{-(l-2\e)n}$ satisfy that $r_n>r_{n+1}=e^{-(l-2\e)}r_n$ for all $n\ge 0$,
we have 
\[
\limsup_{r \to 0}\frac{\log \n_\pi\(B\(\wb_\infty(\o), r\)\)}{\log r} \le \frac{h(\pi)+\e}{l-2\e}
\quad \text{for $\Pb$-almost every $\o \in A$}.
\]
Since $A=A_{\e, N_\e}$ and $\Pb(A_{\e, N_\e}) \ge 1-\e$ for all $\e>0$,
we obtain
\[
\limsup_{r \to 0}\frac{\log \n_\pi\(B\(\xib, r\)\)}{\log r}\le \frac{h(\pi)}{l} \quad \text{for $\n_\pi$-almost every $\xib \in (\partial \G)^2$},
\]
as required.
\qed

Next we show the dimension lower bound.
We use the following lemma.

\begin{lemma}\label{Lem:F}
For every $\e>0$ there exists a Borel set $F_\e$ in $(\partial \G)^2$ such that $\n_\pi(F_\e) \ge 1-\e$ and
\[
\liminf_{r \to 0}\frac{\log \n_\pi\(F_\e \cap B(\xib, r)\)}{\log r} \ge \frac{h(\pi)}{l}-\e \quad \text{for $\n_\pi$-almost every $\xib$ in $(\partial \G)^2$}.
\]
\end{lemma}

\proof
Let $\{\Pb^{\wb_\infty(\o)}\}_{\o \in \O}$ be the conditional probability measures associated with $\s(\wb_\infty)$,
where we have
\[
\Pb=\int_\O \Pb^{\wb_\infty(\o)}\,d\Pb(\o)=\int_{(\partial \G)^2}\Pb^\xib\,d\n_\pi(\xib).
\]
For all $\e>0$ and all positive integer $N$, if we define
\[
A_{\e, N}:=\bigcap_{n \ge N}\Big\{\o \in \O \ : \ d_\times(\wb_n(\o), \Pi_\times(\wb_\infty(\o))(ln))\le \e n, \ \pi_n(\wb_n(\o)) \le e^{-n(h(\pi)-\e)}\Big\},
\]
then by \eqref{Eq:ray2} and \eqref{Eq:Shannon}, there exists an $N_\e$ such that
$\Pb(A_{\e, N_\e}) \ge 1-\e$.
Letting $A:=A_{\e, N_\e}$,
we define
\[
F_\e:=\big\{\xib \in (\partial \G)^2 \ : \ \Pb^\xib(A) \ge \e \big\}.
\]
Since 
\begin{align*}
1-\e \le \Pb(A)=\int_{(\partial \G)^2}\Pb^\xib(A)\,d\n_\pi(\xib)
&=\int_{F_\e}\Pb^\xib(A)\,d\n_\pi(\xib)+\int_{F_\e^{\sf c}}\Pb^\xib(A)\,d\n_\pi(\xib)\\
&\le \n_\pi(F_\e)+\e \n_\pi(F_\e^{\sf c}),
\end{align*}
we have $\n_\pi(F_\e) \ge 1-2\e$.

Let $\zb_n=(z_n^{(1)}, z_n^{(2)})$ be any sequence with $|z_n^{(i)}|=\lfloor ln\rfloor$ for $i=1, 2$.
Note that for $\Pb$-almost every $\y \in A$ and for all $n \ge N_\e$,
if $\wb_\infty(\y) \in \Oc(z_n^{(1)}, R)\times \Oc(z_n^{(2)}, R)$,
then
\[
d_\times(\Pi_\times(\wb_\infty(\y))(ln), \zb_n) \le 2R+C',
\]
for a positive constant $C'$ depending only on the hyperbolicity constant of the metric in $\G$,
and thus 
\[
\y \in A \quad \text{and} \quad \wb_\infty(\y) \in \Oc(z_n^{(1)}, R)\times \Oc(z_n^{(2)}, R) \implies \wb_n(\y) \in B(\zb_n, \e n+C)
\]
where $B(\zb_n, T)=B(z_n^{(1)}, T)\times B(z_n^{(2)}, T)$ for $T \ge 0$
and $C=2R+C'$ for a fixed $R>0$.
This shows that for all $n \ge N_\e$,
\begin{align*}
&\Pb\big(\wb_\infty(\y) \in F_\e \cap \Oc(z_n^{(1)}, R)\times \Oc(z_n^{(2)}, R)\big)\\
&\le \Pb\big(A\cap \{\wb_n(\y) \in B(\zb_n, \e n +C)\}\big)
+\Pb\big(A^{\sf c}\cap \{\wb_\infty(\y) \in F_\e \cap \Oc(z_n^{(1)}, R)\times \Oc(z_n^{(2)}, R)\}\big).
\end{align*}
The second term is estimated as follows:
\begin{align*}
&\Pb\big(A^{\sf c}\cap \{\wb_\infty(\y) \in F_\e \cap \Oc(z_n^{(1)}, R)\times \Oc(z_n^{(2)}, R)\}\big)\\
&=\int_{F_\e \cap \Oc(z_n^{(1)}, R)\times \Oc(z_n^{(2)}, R)}
\Pb^\xib(A^{\sf c})\,d\n_\pi(\xib)
\le (1-\e)\n_\pi\(F_\e \cap \Oc(z_n^{(1)}, R) \times \Oc(z_n^{(2)}, R)\).
\end{align*}
Therefore we obtain
\begin{equation}\label{Eq:LemF1}
\e \n_\pi\big(F_\e \cap \Oc(z_n^{(1)}, R)\times \Oc(z_n^{(2)}, R)\big)
\le \Pb\big(A\cap \{\wb_n(\y) \in B(\zb_n, \e n +C)\}\big),
\end{equation}
for all $n \ge N_\e$.
Moreover, since
\begin{align*}
\Pb\big(A\cap \{\wb_n(\y) \in B(\zb_n, \e n +C)\}\big)
&\le \Pb\big(\pi_n(\wb_n(\y)) \le e^{-n(h(\pi)-\e)}, \ \wb_n(\y) \in B(\zb_n, \e n+C)\big)\\
&\le \#B(\zb_n, \e n+C)\cdot e^{-n(h(\pi)-\e)},
\end{align*}
we have for all $n \ge N_\e$,
\begin{equation}\label{Eq:LemF2}
\Pb\big(A\cap \{\wb_n(\y) \in B(\zb_n, \e n +C)\}\big) \le e^{2D(\e n+C)}e^{-n(h(\pi)-\e)}
\end{equation}
where $D$ is a constant greater than the exponential growth rate of $(\G, d)$, i.e.,
\[
\# B(z_n^{(i)}, T)\le e^{D T} \quad \text{for $i=1, 2$ and for all large enough $T$}.
\]
Combining \eqref{Eq:LemF1} and \eqref{Eq:LemF2}, we obtain
\[
\liminf_{n \to \infty}\frac{\log \n_\pi\(F_\e \cap \Oc(z_n^{(1)}, R)\times \Oc(z_n^{(2)}, R)\)}{-n}
\ge h(\pi)-\e-2D \e.
\]

Let us define $\zb_n:=\Pi_\times(\wb_\infty(\o))(ln)$ for $\Pb$-almost every $\o \in A_{\d, N_\d}$ for all $\d>0$.
By \eqref{Eq:shadow_ball} and \eqref{Eq:prod_ball}, as in a similar way in the last part in the proof of Lemma \ref{Lem:upperbound}, for $\Pb$-almost every $\o \in A_{\d, N_\d}$, 
\[
\liminf_{r \to 0}\frac{\log \n_\pi\(F_\e \cap B(\wb_\infty(\o), r)\)}{\log r} \ge \frac{h(\pi)-\e -2D\e}{l}=\frac{h(\pi)}{l}-C'\e,
\]
for a constant $C'>0$.
Since $\Pb(A_{\d, N_\d}) \ge 1-\d$ for all $\d>0$, we have
\[
\liminf_{r \to 0}\frac{\log \n_\pi\(F_\e\cap B(\xib, r\))}{\log r} \ge \frac{h(\pi)}{l}-C'\e \quad \text{for $\n_\pi$-almost every $\xib$ in $(\partial \G)^2$}.
\]
Replacing $\e$ by a small enough constant yields the claim as stated.
\qed

\begin{lemma}\label{Lem:lowerbound}
For $\n_\pi$-almost all $\xib$ in $(\partial \G)^2$,
\[
\liminf_{r \to 0}\frac{\log \n_\pi\(B(\xib, r)\)}{\log r} \ge \frac{h(\pi)}{l}.
\]
\end{lemma}

\proof
For all $\e>0$, let $F:=F_\e$ be the Borel set in Lemma \ref{Lem:F}.
For the hyperbolic metric space $(\G, d)$ and the boundary $(\partial \G, q)$,
we have that for each $0<\a<1$, 
the space $(\partial \G, q^\a)$ admits a bi-Lipschitz embedding into the Euclidean space $\R^n$ for some $n$ \cite[Theorem 9.2]{BonkSchramm}.
Hence there exists a bi-Lipschitz map $\f=(\f_1, \f_2)$,
\[
\f: ((\partial \G)^2, q_\times^\a) \to \R^n \times \R^n, \quad (\x^{(1)}, \x^{(2)}) \mapsto (\f_1(\x^{(1)}), \f_2(\x^{(2)})),
\]
i.e., for a constant $L>0$,
\[
\frac{1}{L}q_\times(\xib_1, \xib_2)^\a\le \|\f(\xib_1)-\f(\xib_2)\|_{\R^{2n}}\le Lq_\times(\xib_1, \xib_2)^\a
\]
for all $\xib_i \in (\partial \G)^2$, $i=1, 2$,
where $\|\cdot\,\|_{\R^{2n}}$ denotes the standard Euclidean norm in $\R^{2n}=\R^n \times \R^n$.
By the Lebesgue density theorem for the Borel measure $\f_\ast \n_\pi$ in $\R^{2n}$,
we have
\[
\lim_{r \to 0}\frac{\f_\ast\n_\pi\(\f(F)\cap B_{\R^{2n}}(\f(\xib), r)\)}{\f_\ast \n_\pi\(B_{\R^{2n}}(\f(\xib), r)\)}=1 \quad \text{for $\n_\pi$-almost every $\xib \in F$},
\]
where $B_{\R^{2n}}(x, r)$ stands for the standard Euclidean (closed) ball in $\R^{2n}$.
This implies that
\[
\liminf_{r \to 0}\frac{\n_\pi\(F \cap B(\xib, (L r)^{1/\a})\)}{\n_\pi\(B(\xib, (r/L)^{1/\a})\)} \ge 1 \quad \text{for $\n_\pi$-almost all $\xib \in F$},
\]
and for $\n_\pi$-almost all $\xib \in (\partial \G)^2$, there exist positive constants $c(\xib)>0$ and $r(\xib)>0$ such that
\[
\n_\pi\(F\cap B(\xib, L^{2/\a}r)\) \ge c(\xib)\n_\pi\(B(\xib, r)\) \quad \text{for all $0<r<r(\xib)$}.
\]
Therefore we obtain
\[
\liminf_{r \to 0}\frac{\log \n_\pi\(B(\xib, r)\)}{\log r}\ge \liminf_{r \to 0}\frac{\log \n_\pi\(F\cap B(\xib, r)\)}{\log r}
\quad \text{for $\n_\pi$-almost all $\xib \in F$}.
\]
Lemma \ref{Lem:F} implies that
\[
\liminf_{r \to 0}\frac{\log \n_\pi\(B(\xib, r)\)}{\log r} \ge \frac{h(\pi)}{l}-\e
\quad \text{for $\n_\pi$-almost every $\xib \in F$},
\]
and since $F=F_\e$ and $\n_\pi(F_\e) \ge 1-\e$ for all $\e>0$,
\[
\liminf_{r \to 0}\frac{\log\n_\pi\(B(\xib, r)\)}{\log r} \ge \frac{h(\pi)}{l} \quad \text{for $\n_\pi$-almost every $\xib$ in $(\partial \G)^2$},
\]
concluding the claim.
\qed

\proof[Proof of Theorem \ref{Thm:dim}]
Lemmas \ref{Lem:upperbound} and \ref{Lem:lowerbound} show the claim.
\qed

\section{Continuity of entropy}\label{Sec:conti}

For a countable group $\G$ (in particular we discuss a product of word hyperbolic groups),
we endow the set of probability measures on $\G$ with the topology induced by the total variation distance.
Note that for all probability measure $\m$ and all sequence of probability measures $\{\m_{(i)}\}_{i=0}^\infty$ 
we have
\[
\|\m_{(i)}-\m\|_\TV \to 0 \quad \text{as $i \to \infty$},
\]
if and only if
$\m_{(i)}(x) \to \m(x)$ as $i \to \infty$
for each $x \in \G$.
Fix a left-invariant metric $d$ on $\G$ with finite exponential growth rate and let $|x|=d(o, x)$ for $x \in \G$.
For a finite set $K$ in $\G$,
let
\[
\Eb_\m [|x|: \G \setminus K]:=\sum_{x \in \G \setminus K}|x|\,\m(x)
\]

Erschler and Kaimanovich have shown that the continuity of $h(\m)$ in $\m \in \Mcc$ under some general conditions \cite{ErschlerKaimanovich}.
We say that a set $\Mcc$ of probability measures on $\G$ satisfies \textit{uniform first moment condition} if 
\begin{equation}\tag{M}\label{M}
\sup_{\m \in \Mcc}\Eb_\m[|x|:\G \setminus K_n] \to 0 \quad \text{as $n \to \infty$},
\end{equation}
for all sequence of finite sets $\{K_n\}_{n=0}^\infty$ with $\bigcup_{n=0}^\infty K_n=\G$.
We assume that there exists a pair of Borel $\G$-spaces $B$, $\check B$ such that the $\G$-space $\check B \times B$ with the diagonal action admits a $\G$-invariant Borel set $\L$ in $\check B \times B$ and a $\G$-equivariant map $S$ assigning to $(\check \x, \x) \in \L$ a proper subset (\textit{strip}) in $\G$. 
Let us say that the strips $S(\check \x, \x)$ given by the map $S$ satisfy \textit{uniform subexponential growth} if
\begin{equation}\tag{G}\label{G}
\sup_{(\check \x, \x)\in \L}\frac{1}{n}\log \#\(B(o, n) \cap S(\check \x, \x)\) \to 0 \quad \text{as $n \to \infty$}.
\end{equation}
For non-negative real $R$, letting $S_R(\check \x, \x)$ be the $R$-neighborhood of $S(\check \x, \x)$,
we define
\[
\L_R:=\big\{(\check \x, \x) \in \L \ : \ o \in S_R(\check \x, \x)\big\}.
\]
Note that the union of $\L_R$ over $R$ covers $\L$.
If a pair of Borel $\G$-spaces $B$ and $\check B$ admits a probability measure $\n_\m$ on $B$ (resp.\ $\n_{\check \m}$ on $\check B$) for which $(B, \n_\m)$ (resp.\ $(\check B, \n_{\check \m})$) is a $\m$- (resp.\ $\check \m$-) boundary (where $\check \m(x):=\m(x^{-1})$ for $x \in \G$), and further
\begin{equation}\tag{S}\label{S}
\inf_{\m \in \Mcc}\n_{\check \m} \times \n_\m(\L_R) \to 1 \quad \text{as $R \to \infty$},
\end{equation}
then we say that $\Mcc$ satisfies \textit{uniform strip condition}.

\begin{theorem}[Theorem 1 in \cite{ErschlerKaimanovich}]\label{Thm:EK}
If a set $\Mcc$ of probability measures on $\G$ and a map $S$ satisfy the conditions \eqref{M}, \eqref{G} and \eqref{S},
then the function $\Mcc \to \R$, $\m \mapsto h(\m)$ is continuous.
\end{theorem}

Theorem \ref{Thm:EK} applies to word hyperbolic groups and their products with sequences of probability measures:
In the case when $\G$ is a word hyperbolic group,
for a sequence of probability measures $\{\m_{(i)}\}_{i=0}^\infty$ on $\G$
with uniform first moment converging to a probability measure $\m$,
we have $h(\m_{(i)}) \to h(\m)$ as $i \to \infty$ \cite[Theorem 2]{ErschlerKaimanovich}.
Actually, it suffices to consider the case when $\G$ is a non-elementary word hyperbolic group and the limiting probability measure $\m$ is non-elementary.
One may take $B=\check B=\partial \G$ the Gromov boundary endowed with the harmonic measures $\n_\m$ and $\n_{\check \m}$, respectively, and
$\L:=(\partial \G)^2 \setminus\{{\rm diagonal}\}$, which is open in the product $(\partial \G)^2$.
Furthermore, for $(\check \x, \x) \in \L$ the strip $S(\check \x, \x)$ is defined as the union of bi-infinite geodesics connecting $\check \x$ and $\x$ in a Cayley graph of $\G$.
The condition \eqref{G} is satisfied since 
\[
\#(B(o, n)\cap S(\check \x, \x))=O(n),
\]
where the implied constant depends only on the Cayley graph. 
Furthermore the condition \eqref{S} is satisfied.
Indeed, the harmonic measure $\n_\m$ is the unique $\m$-stationary measure on $\partial \G$ and the measures $\n_{\m_{(i)}}$ weakly converge to $\n_\m$ as $i \to \infty$, and $\n_\m$ is supported on an open set $\L$.
Letting $\Lambda_R^\circ$ denote the interior of $\Lambda_R$, we have
\[
\liminf_{k \to \infty}\check \nu_{i_k}\times \nu_{i_k}(\Lambda_R^\circ) \ge \nu_{\check \mu} \times \nu_\mu(\Lambda_R^\circ),
\]
for every subsequence $\check \nu_{i_k}\times \nu_{i_k}$ of $\nu_{\check \mu_{(i)}} \times \nu_{\mu_{(i)}}$.
Noting that $\Lambda_R^\circ$ increases and exhausts $\Lambda$ as $R$ grows,
we have \eqref{S} (cf.\ \cite[Lemma 3]{ErschlerKaimanovich}).

In the case when the group is a product $\G \times \G$ of word hyperbolic groups and a probability measure $\pi$,
one may take $B=\check B=(\partial \G)^2$ and 
\[
\L=\big\{(\check \xib, \xib) \in \check B\times B \ : \ \check \x^{(i)}\neq \x^{(i)} \ \text{for $i=1, 2$}\big\},
\]
where $\check \xib=(\check \x^{(1)}, \check\x^{(2)})$ and $\xib=(\x^{(1)}, \x^{(2)})$, and $\L$ is open in $\check B \times B$.
The strip is defined by
\[
S(\check \xib, \xib)=S(\check \x^{(1)}, \x^{(1)})\times S(\check \x^{(2)}, \x^{(2)}),
\]
and we have
\[
\# \(B(o, n)\cap S(\check \xib, \xib)\)=O(n^2).
\]
This shows that \eqref{G} holds.
Moreover we have \eqref{S} for a sequence of probability measures $\pi_{(i)}$ on $\Gamma \times \Gamma$ 
since $\n_\pi$ is the unique $\pi$-stationary measure on $(\partial \G)^2$ and supported on an open set $\L$ as in the case on $\Gamma$ presented above.

\begin{corollary}\label{Cor:conti}
For a word hyperbolic group $\G$ and a non-elementary probability measure $\m$ with finite first moment,
the asymptotic entropy $h(\pi^\rho)$ is continuous in $\rho \in [0, 1]$.
\end{corollary}

\proof
The set of probability measures $\pi^\rho=\rho \m\times \m+(1-\rho)\m_{\diag}$ on $\G \times \G$ has uniform finite first moment if $\m$ has finite first moment.
By the discussion above, the conditions \eqref{M}, \eqref{G} and \eqref{S} are satisfied for $\{\pi^{\rho_i}\}_{i=0}^\infty$ with every sequence $\{\rho_i\}_{i=0}^\infty$ converging to $\rho$ in $[0, 1]$ as $i \to \infty$, and thus Theorem \ref{Thm:EK} implies the claim.
\qed

\section{Proofs of Theorems \ref{Thm:1} and \ref{Thm:2}}\label{Sec:proofs}

\proof[Proof of Theorem \ref{Thm:1}]
For probability measures $\n_1$ and $\n_2$ on $(\G \cup \partial \G)^2$,
the total variation distance is defined by
\[
\|\n_1-\n_2\|_\TV:=\sup \big\{|\n_1(A)-\n_2(A)| \ : \ \text{$A$ is Borel in $(\G \cup \partial \G)^2$} \big\}.
\]
For each $0\le \rho \le 1$, a $\pi^\rho$-random walk $\{\wb_n\}_{n=0}^\infty$ converges to $\wb_\infty$ in $(\partial \G)^2$ as $n \to \infty$ in $(\G \cup \partial \G)^2$, $\Pb$-almost surely, and
the distribution $\pi^\rho_n$ converges weakly to the harmonic measure $\n_{\pi^\rho}$ (see the beginning of Section \ref{Sec:dim}).
Therefore for $0 \le \rho \le 1$, we have
\begin{equation}\label{Eq:Thm1liminf}
\liminf_{n \to \infty}\|\pi^\rho_n -\m_n\times \m_n\|_{\TV} \ge \|\n_{\pi^\rho}-\n_\m\times \n_\m\|_{\TV}.
\end{equation}
For all $0\le \rho<1$, we have $\n_{\pi^\rho}\neq \n_\m\times \n_\m$.
Indeed, suppose that $\n_{\pi^\rho}=\n_\m\times \n_\m$ for some $0\le \rho<1$, then we have
\[
\pi^\rho\ast (\n_\m\times\n_\m)=\n_\m\times \n_\m
\]
since $\n_{\pi^\rho}$ is the $\pi^\rho$-stationary measure on $(\partial \G)^2$ (cf.\ \eqref{Eq:pi-stationary} and \eqref{Eq:mu-stationary} in Section \ref{Sec:dim}).
Noting that $\n_\m$ is the $\m$-stationary measure on $\partial \G$, we have
\[
\rho (\n_\m \times \n_\m)+(1-\rho)\sum_{x \in \G}\m_{\diag}(x)(x\n_\m \times x \n_\m)=\n_\m \times \n_\m,
\]
and
\[
\m_{\diag}\ast (\n_\m\times \n_\m)=\n_\m \times \n_\m.
\]
This shows that $\n_\m\times \n_\m$ is the $\m_\diag$-stationary (harmonic) measure by the uniqueness.
However, the harmonic measure for $\m_{\diag}$ is supported on the diagonal in $(\partial \G)^2$
and $\n_\m$ is non-atomic on $\partial \G$ since $\gr(\m)$ is non-elementary, we have $\n_{\pi^\rho}\neq \n_\m \times \n_\m$, yielding a contradiction.
Therefore for all $0\le \rho<1$, we have $\|\n_{\pi^\rho}-\n_\m \times \n_\m\|_\TV>0$,
and thus by \eqref{Eq:Thm1liminf}, 
\[
\liminf_{n \to \infty}\|\pi^\rho_n-\m_n\times \m_n\|_{\TV}>0,
\]
as claimed.
\qed

\proof[Proof of Theorem \ref{Thm:2}]
As in the same way in the beginning of the proof of Theorem \ref{Thm:1},
for all $0 \le \rho, \rho' \le 1$, we have
\[
\liminf_{n \to \infty}\|\pi^\rho_n-\pi^{\rho'}_n\|_{\TV}\ge \|\n_{\pi^\rho}-\n_{\pi^{\rho'}}\|_{\TV}.
\]
Theorem \ref{Thm:dim} shows that for the Borel set
\[
E_\rho:=\left\{\xib \in (\partial \G)^2 \ : \ \lim_{r \to 0}\frac{\log \n_{\pi^\rho}(B(\xib, r))}{\log r}=\frac{h(\pi^\rho)}{l}\right\},
\]
we have $\n_{\pi^\rho}(E_\rho)=1$.
By Corollary \ref{Cor:conti}, the function $\rho\mapsto h(\pi^\rho)$ for $\rho \in [0, 1]$ is continuous.
Furthermore, 
\[
h(\m)=h(\pi^0)\le h(\pi^\rho) \le h(\pi^1)=2h(\m) \quad \text{for all $0 \le \rho \le 1$}, 
\]
and $h(\m)>0$ since $\gr(\m)$ is non-elementary. 
Hence
there exists $0<\rho_\ast\le 1$ such that $h(\pi^\rho)<h(\pi^1)$ for all $0 \le \rho<\rho_\ast$.
This shows that for all $0\le \rho<\rho_\ast$, we have $\n_{\pi^\rho}(E_\rho)=1$ and $(\n_\m\times \n_\m)(E_\rho)=0$, implying that $\n_{\pi^\rho}$ and $\n_\m \times \n_\m$ are mutually singular and 
$\|\n_{\pi^\rho}-\n_\m\times \n_\m\|_{\TV}=1$.
Therefore we have for all $0\le \rho<\rho_\ast$,
\[
\liminf_{n \to \infty}\|\pi^\rho_n-\m_n \times \m_n\|_{\TV}=1,
\]
and $\lim_{n \to \infty}\|\pi^\rho_n-\m_n\times \m_n\|_{\TV}=1$, as required.
\qed

\begin{remark}\label{Rem:entropy}
The proof of Theorem \ref{Thm:2} shows that if $h(\pi^\rho)\neq h(\pi^{\rho'})$ for $0 \le \rho, \rho' \le 1$,
then 
\[
\|\pi_n^\rho-\pi_n^{\rho'}\|_{\TV} \to 1 \quad \text{as $n \to \infty$},
\]
and $\nu_{\pi^\rho}$ and $\nu_{\pi^{\rho'}}$ are mutually singular.
\end{remark}

\appendix

\section{A proof of Theorem \ref{Thm:EK}}\label{Appendix}

In this section, $\Eb_\n$ denotes the expectation for a probability measure $\n$.
Let $\G$ be a countable group endowed with a probability measure $\m$ of finite entropy, i.e., $H(\m)<\infty$.
For the $\m$-random walk $\{w_n\}_{n=0}^\infty$ starting from $o$ on $\G$,
let us consider the probability measure space $(\G^{\Z_+}, \Fc, \Pb)$, where $\Fc$ is the $\s$-field generated by the cylinder sets and $\Pb$ is the distribution of $\{w_n\}_{n=0}^\infty$.

For each positive integer $n$, let $\a_1^n$ be the measurable partition on $\G^{\Z_+}$ where $\ob=(\o_i)_{i=0}^\infty$ and $\ob'=(\o_i')_{i=0}^\infty$ belong to the same set if and only if $\o_i=\o_i'$ for all $0 \le i \le n$.
For any sub $\s$-field $\Ac$ in $\Fc$,
the conditional entropy is defined by
\begin{equation*}
H_\Pb(\a_1^n \mid \Ac):=\Eb_\Pb\Big[-\sum_{B \in \a_1^n}\Pb(B \mid \Ac)\log \Pb(B\mid \Ac)\Big],
\end{equation*}
where $\Pb(\,\cdot \mid \Ac)$ stands for the conditional probability measure with respect to $\Ac$.
The \textit{tail} $\s$-field is defined by $\Tc:=\bigcap_{n=0}^{\infty}\s(w_n, w_{n+1}, \dots)$.
In the case when $\Ac=\Tc$, letting $\a:=\a_1^1$, we have
\begin{equation}\label{Eq:KVT}
H_\Pb(\a \mid \Tc)=H(\m)-h(\m),
\end{equation}
\cite[cf.\ Proof of Theorem 1.1]{KaimanovichVershik}.

The group $\G$ acts on $\G^{\Z_+}$ by $x(\o_n)_{n=0}^\infty=(x\o_n)_{n=0}^\infty$ for $x \in \G$.
The \textit{stationary $\s$-field} $\Sc$ is the sub $\s$-field of $\Fc$ generated by shift-invariant measurable sets, where the shift is defined by
$(\o_n)_{n=0}^\infty \mapsto (\o_{n+1})_{n=0}^\infty$ on $\G^{\Z_+}$.
Note that $\Sc$ is $\G$-invariant, i.e., if $A \in \Sc$, then $x A \in \Sc$ for all $x \in \G$.
By definition, we have $\Sc \subset \Tc$, and it is known that their $\Pb$-completions coincide, i.e., $\Sc=\Tc$ mod $\Pb$ \cite[Section 7.0]{KaimanovichVershik} (where it is crucial that the initial state $w_0$ is a point).
Therefore by \eqref{Eq:KVT},
\begin{equation}\label{Eq:KV}
H_{\Pb}(\a \mid \Sc)=H(\m)-h(\m).
\end{equation}

For each $\G$-invariant sub $\s$-field $\Ac$ of $\Sc$, let $\Pb^\x(\,\cdot\,)=\Pb(\,\cdot\mid \Ac)(\x)$ for $\Pb$-almost every $\x \in \G^{\Z_+}$, and
let
\[
\m^\x_n(x):=\Pb^\x(w_n=x) \qquad \text{for $x \in \G$ and $n \ge 0$},
\]
and we define the entropy of conditional process: For $n \ge 0$, let
\begin{equation}\label{Eq:cond_ent}
H(\m_n^\x):=\Eb_\Pb\Big[-\sum_{x \in \G}\m_n^\x(x)\log \m_n^\x(x)\Big].
\end{equation}
This yields
$H(\m_n^\x)=H(\m_n)+n(H(\a\mid \Ac)-H(\m))$,
and in particular, in the case when $\Ac=\Sc$, by \eqref{Eq:KV},
we obtain for all $n \ge 0$,
\begin{equation}\label{Eq:conditional_entropy}
H(\m_n^\x)=H(\m_n)-n h(\m),
\end{equation}
\cite[Sections 3 and 4]{Kaimanovich-hyperbolic}.

We write $B_R=B(o, R)$ for simplicity of notations.

\begin{lemma}\label{Lem:conti}
For a set of probability measures $\Mcc$ on $\G$ with uniform first moment condition,
the function $\Mcc \to \R$, 
\[
\m \mapsto H(\m)
\]
is continuous.
\end{lemma}

\proof
For the exponential growth rate $v(\G, d)$ for $(\G, d)$,
let us fix $D>v(\G, d)$ and define
\[
A:=\big\{x \in \G \ : \ \m(x) \ge e^{-D|x|}\big\}.
\]
For the ball $K=B_N$ in $\G$ with $N$ large enough,
decomposing the sum
\[
H(\m)=-\sum_{x \in K}\m(x)\log \m(x)
-\sum_{x \in A \cap K^{\sf c}}\m(x)\log \m(x)-\sum_{x \in A^{\sf c} \cap K^{\sf c}}\m(x)\log \m(x),
\]
we estimate the second and third terms.
First, we have
\begin{align*}
-\sum_{x \in A \cap K^{\sf c}}\m(x)\log \m(x) 
&\le -\sum_{x \in K^{\sf c}}\m(x)\log e^{-D|x|} \le D\sup_{\m \in \Mcc}\Eb_\m[|x| : \G \setminus K].
\end{align*}
Second, letting $S_k:=\{x \in \G \ : \ k< |x|\le k+1\}$ for non-negative integers $k$,
\begin{align*}
-\sum_{x \in A^{\sf c} \cap K^{\sf c}}\m(x)\log \m(x)
&=-\sum_{k=N}^\infty \sum_{x \in S_k \cap A^{\sf c}}\m(x)\log \m(x)
\le \sum_{k=N}^\infty \#S_k \cdot D(k+1) e^{-Dk} \le C e^{-D' N}
\end{align*}
for constants $0<D'<D-v(\G, d)$ and $C>0$, for all large enough $N$,
where in the first inequality we have used $-\m(x)\log \m(x) \le -e^{-D|x|}\log e^{-D|x|}$ for $x \in A^{\sf c}$ and $|x|$ large enough.
Finally, we obtain
\[
\sup_{\m \in \Mcc}\left|H(\m)+\sum_{x \in B_N}\m(x)\log \m(x)\right|
\le D\sup_{\m \in \Mcc}\Eb_\m\Big[|x| : \G \setminus B_N\Big]+C e^{-D' N}.
\]
This shows that
\[
-\sum_{x \in B_N}\m(x)\log \m(x) \to H(\m) \quad \text{uniformly on $\m \in \Mcc$ as $N \to \infty$},
\]
implying that $\m \mapsto H(\m)$ is continuous on $\Mcc$.
\qed

\begin{lemma}\label{Lem:entropy_outside}
In the same setting as in Lemma \ref{Lem:conti},
for all $L>4$ and for all positive integer $n$,
\begin{equation}\label{Eq:Lem:entropy1}
\sup_{\m \in \Mcc}\Eb_{\m_n}\Big[|x|:\G \setminus B_{nL}\Big] 
\le n \sup_{\m \in \Mcc}\Eb_\m \Big[|x|: \G \setminus B_{\sqrt{L}}\Big]
+\frac{2n}{\sqrt{L}}\sup_{\m \in \Mcc}\Eb_\m|x|,
\end{equation}
and
\begin{equation}\label{Eq:Lem:entropy2}
\sup_{\m \in \Mcc}\Big(-\sum_{x \in \G \setminus B_{nL}}\m_n(x)\log \m_n(x)\Big)
\le D\sup_{\m \in \Mcc}\Eb_{\m_n}\Big[|x| : \G \setminus B_{nL}\Big]+C e^{-D' nL},
\end{equation}
where $C$, $D$ and $D'$ are positive constants independent of $n$ and $L$.
Moreover, for each $n>0$, the function $\Mcc \to \R$, $\m \mapsto H(\m_n)$ is continuous.
\end{lemma}

\proof
We use the same notation as in the proof of Lemma \ref{Lem:conti} and obtain \eqref{Eq:Lem:entropy2} in the same way 
for each positive integer $n$ and for all $L>4$.
Let us show \eqref{Eq:Lem:entropy1}.
Note that 
\begin{align*}
\Eb_{\m_n}\Big[|x|:\G \setminus B_{nL}\Big]
=\Eb_\Pb \Big[|w_n|\1_{\{|w_n|>nL\}}\Big]
&\le \Eb_\Pb\Big[\Big(\sum_{i=1}^n|x_i|\Big)\1_{\{\sum_{i=1}^n |x_i|>nL\}}\Big]\\
&=n \Eb_\Pb \Big[|x_1|\1_{\{\sum_{i=1}^n |x_i|>nL\}}\Big],
\end{align*}
where the last equality follows 
since $x_1, \dots, x_n$ are independent and identically distributed.
Moreover, we have
\begin{align*}
\Eb_\Pb\Big[|x_1|\1_{\{\sum_{i=1}^n |x_i|>nL\}}\Big]
=\Eb_\Pb\Big[|x_1|\1_{\{|x_1|>\sqrt{L}, \ \sum_{i=1}^n |x_i|>nL\}}\Big]
+\Eb_\Pb\Big[|x_1|\1_{\{|x_1|\le \sqrt{L}, \ \sum_{i=1}^n |x_i|>nL\}}\Big],
\end{align*}
where the first term is at most
\[ 
\Eb_\m\Big[|x_1|\1_{\{|x_1|>\sqrt{L}\}}\Big]
\le \sup_{\m \in \Mcc}\Eb_\m \Big[|x|: \G \setminus B_{\sqrt{L}}\Big],
\]
and the second term is at most
\[
\sqrt{L}\Pb\Big(\sum_{i=2}^n |x_i|>nL-\sqrt{L}\Big)\le \frac{\sqrt{L}(n-1)}{nL-\sqrt{L}}\Eb_\m|x|
\le \frac{\sqrt{L}}{L-\sqrt{L}/n}\sup_{\m \in \Mcc}\Eb_\m|x|,
\]
by the Markov inequality.
Therefore for all $L>4$ and for all $n>0$, we have $\sqrt{L}/n<L/2$ and
\[
\sup_{\m \in \Mcc}\Eb_{\m_n}\Big[|x|:\G \setminus B_{nL}\Big] 
\le n \sup_{\m \in \Mcc}\Eb_\m \Big[|x|: \G \setminus B_{\sqrt{L}}\Big]
+\frac{2n}{\sqrt{L}}\sup_{\m \in \Mcc}\Eb_\m|x|,
\]
yielding \eqref{Eq:Lem:entropy1}.
Finally, since $\|\m_n-\m_n'\|_\TV\le n\|\m-\m'\|_\TV$ for $\m$, $\m' \in \Mcc$, the term 
\[
-\sum_{x \in B_{nL}}\m_n(x)\log \m_n(x)
\] 
is continuous in $\m$ for each fixed $n$,
and 
converges to $H(\m_n)$ uniformly on $\m$ in $\Mcc$ as $L \to \infty$ by \eqref{Eq:Lem:entropy2},
the last statement holds.
\qed

Recall the notations from \cite[Section 6]{Kaimanovich-hyperbolic}:
Let $(\G^\Z, \wbar \Pb)$ be the probability measure space of bilateral paths $\wbar \ob=(\o_i)_{i \in \Z}$ with $\omega_0=o$.
The space is identified with the product space via the map $\wbar \ob \mapsto (\check \ob, \ob)$ from $(\G^\Z, \wbar \Pb)$ to $(\G^{\Z_+}, \check \Pb)\times (\G^{\Z_+}, \Pb)$ where $\check \ob=(\o_{-i})_{i \in \Z_+}$ and $\check \Pb$ is the distribution of $\check \m$-random walk starting from $o$.
We denote by $\wbar U$ the probability measure preserving transformation on $(\G^\Z, \wbar \Pb)$ induced from the Bernoulli shift in the space of increments, more explicitly,
\[
(\wbar U^k \ob)_n=\o_k^{-1}\o_{n+k} \quad \text{for $\ob=(\o_i)_{i \in \Z} \in \G^\Z$ and for $k, n \in \Z$}.
\]
Given $\G$-equivariant measurable maps ${\rm bnd_+}:\G^{\Z_+} \to B$ and ${\rm bnd_-}:\G^{\Z_+} \to \check B$ for the $\m$-boundary $B$ and for the $\check \m$-boundary $\check B$,
let us define $\Pi_+:\G^\Z \to B$ by $\wbar \ob=(\check \ob, \ob) \mapsto {\rm bnd_+}(\ob)$
and $\Pi_-:\G^\Z \to \check B$ by $\wbar \ob=(\check \ob, \ob) \mapsto {\rm bnd_-}(\check \ob)$.
Note that $\n_\m={\Pi_+}_\ast\wbar \Pb$ and $\n_{\check \m}={\Pi_-}_\ast\wbar \Pb$.

\proof[Proof of Theorem \ref{Thm:EK}]
The condition \eqref{S} implies that 
\[
\e_R:=1-\inf_{\m \in \Mcc}\n_{\check \m}\times \n_\m(\L_R) \to 0 \quad \text{as $R \to \infty$},
\]
where 
$\L_R=\big\{(\check \x, \x) \in \L \ : \ o \in S_R(\check \x, \x)\big\}$ for $R \ge 0$,
and thus
\begin{align*}
\wbar \Pb\(o \in S_R(\Pi_- \wbar \ob, \Pi_+\wbar \ob)\)=
\n_{\check \m}\times \n_\m\((\check \x, \x) \in \L \ : \ o \in S_R(\check \x, \x)\)=\n_{\check \m}\times \n_{\m}\(\L_R\) \ge 1-\e_R,
\end{align*}
uniformly on $\m \in \Mcc$.
Moreover, since the map $S$ is $\G$-equivariant and $\wbar \Pb$ is $\wbar U$-invariant,
we have
\begin{align*}
&\wbar \Pb\(\o_n \in S_R(\Pi_-\wbar \ob, \Pi_+\wbar \ob)\)
=\wbar \Pb\(o \in \o_n^{-1}S_R(\Pi_-\wbar \ob, \Pi_+\wbar \ob)\)\\
&\qquad \qquad \qquad=\wbar \Pb\(o \in S_R(\Pi_- \wbar U^n \wbar \ob, \Pi_+ \wbar U^n \wbar \ob)\) =\wbar \Pb\(o \in S_R(\Pi_- \wbar \ob, \Pi_+ \wbar \ob)\) \ge 1-\e_R.
\end{align*}
Therefore, disintegrating the measure,
\[
\wbar \Pb\(\o_n \in S_R(\Pi_- \wbar \ob, \Pi_+ \wbar \ob)\)
=\int_{\L}\Pb^{\x}\(\o_n \in S_R(\check \x, \x)\)\,d\n_{\check \m}d\n_\m
=\int_\L \m_n^\x(S_R(\check \x, \x))d\n_{\check \m}d\n_\m,
\]
we obtain
\begin{equation}\label{Eq:SR}
\int_\L \m_n^\x(S_R(\check \x, \x))d\n_{\check \m}d\n_\m \ge 1-\e_R.
\end{equation}
By Lemma \ref{Lem:entropy_outside} \eqref{Eq:Lem:entropy1}, for all $L>4$ and for all $n>0$,
\[
\sup_{\m \in \Mcc}\Eb_{\m_n}\Big[|x|:\G \setminus B_{nL}\Big] 
\le n \sup_{\m \in \Mcc}\Eb_\m \Big[|x|: \G \setminus B_{\sqrt{L}}\Big]
+\frac{2n}{\sqrt{L}}\sup_{\m \in \Mcc}\Eb_\m|x|,
\]
and thus letting
\[
\e_L:=\sup_{\m \in \Mcc}\Eb_\m \Big[|x|: \G \setminus B_{\sqrt{L}}\Big] \quad \text{and} \quad C_\m:=\sup_{\m \in \Mcc}\Eb_\m|x|,
\]
we have $\e_L \to 0$ as $L \to \infty$ by the condition \eqref{M}, and by the Markov inequality,
\[
\Pb(|w_n| >n L) \le \frac{\Eb_{\m_n}\Big[|x|:\G \setminus B_{nL}\Big]}{n L} 
\le \frac{1}{L}\(\e_L+\frac{2C_\m}{\sqrt{L}} \)=\frac{\tilde \e_L}{L},
\quad \text{where $\tilde \e_L:=\e_L+\frac{2C_\m}{\sqrt{L}}$},
\]
and $\tilde \e_L \to 0$ as $L \to \infty$.
Hence for all $n>0$ and all $\m \in \Mcc$,
we have $\Pb(|w_n| \le nL)\ge 1-\tilde \e_L/L$,
and
this yields by disintegration,
\begin{equation}\label{Eq:BnL}
\int_B \m_n^\x(B_{nL})\,d\n_\m \ge 1-\frac{\tilde \e_L}{L}.
\end{equation}
For all $\e>0$, let us take any $L$ and $R$ satisfying that 
$\tilde \e_L<\e$ and $\e_R<\e/L$.
Noting that
\[
\m_n^\x\(B_{nL}\cap S_R(\check \x, \x)\) \ge \m_n^\x(B_{nL})+\m_n^\x(S_R(\check \x, \x))-1,
\]
we have by \eqref{Eq:SR} and \eqref{Eq:BnL},
\begin{align}\label{Eq:1-2e}
\int_\L \m_n^\x\(B_{nL}\cap S_R(\check \x, \x)\)\,d\n_{\check \m}d\n_\m
&\ge \int_B\m_n^\x\(B_{nL}\)\,d\n_\m+\int_\L \m_n^\x\(S_R(\check \x, \x)\)\,d\n_{\check \m}d\m_\m-1\nonumber\\
&\ge 1-\frac{\tilde \e_L}{L}+1-\e_R-1 \ge 1-\frac{2\e}{L}.
\end{align}

Furthermore, the condition \eqref{G} implies that for all $R>0$ and all $L>0$,
there exists a sequence of positive reals $\f_{n, R, L}$ such that
\begin{equation}\label{Eq:f}
\sup_{(\check \x, \x)\in \L}\#\(S_R(\check \x, \x)\cap B_{nL}\) \le \f_{n, R, L} \quad \text{for all $n>0$},
\end{equation}
and
$(1/n)\log \f_{n, R, L} \to 0$ as $n \to \infty$.

Let us estimate the conditional entropy $H(\m_n^\x)$ in \eqref{Eq:cond_ent}.
For the simplicity of notations,
let
\[
G_1:=S_R(\check \x, \x)\cap B_{nL}, \quad G_2:=B_{nL}\setminus G_1, \quad \text{and} \quad G_3:=\G\setminus B_{nL}.
\]
First we have by the Jensen inequality and by \eqref{Eq:f},
\begin{align*}
-\sum_{x \in G_1}\m_n^\x(x)\log \m_n^\x(x)
&\le \m_n^\x(G_1)\log \#G_1 -\m_n^\x(G_1)\log \m_n^\x(G_1)\\
&\le \log \f_{n, R, L}-\m_n^\x(G_1)\log \m_n^\x(G_1),
\end{align*}
and thus,
\begin{align}\label{Eq:G1}
\int_\L -\sum_{x \in G_1}\m_n^\x(x)\log \m_n^\x(x)\,d\n_{\check \m}d\n_\m
\le \log \f_{n, R, L}+e^{-1},
\end{align}
where we have used $-x \log x \le e^{-1}$ for $0 \le x \le 1$.

Second by \eqref{Eq:1-2e}, we have $\Pb(w_n \in G_2)\le 2\e/L$,
and thus by the Jensen inequality,
\[
-\sum_{x \in G_2}\m_n^\x(x)\log \m_n^\x(x) \le 
\m_n^\x(G_2)\log \# G_2-\m_n^\x(G_2)\log \m_n^\x(G_2),
\]
and as in a similar way to \eqref{Eq:G1},
\begin{align}\label{Eq:G2}
\int_\L-\sum_{x \in G_2}\m_n^\x(x)\log \m_n^\x(x) d\n_{\check \m}d\n_\m
&\le \Pb(w_n \in G_2)\log \#B_{nL}+e^{-1}\nonumber\\
&\le \(\frac{2\e}{L}\)nL D +e^{-1}=2\e n D+e^{-1},
\end{align}
for $D>v(\G, d)$ and all large enough $n$. 
By \eqref{Eq:G1} and \eqref{Eq:G2}, noting that $G_1 \cup G_2=B_{nL}$, we have
\begin{align}\label{Eq:G1G2}
\int_{\L}-\sum_{x \in B_{nL}}\m_n^\x(x)\log \m_n^\x(x)\,d\n_{\check \m}d\n_\m \le \log \f_{n, R, L}+2\e n D+2e^{-1},
\end{align}
for all large enough $n$.

Finally we obtain on $G_3=B_{nL}^{\sf c}$,
\[
\int_{\L}-\sum_{x \in G_3}\m_n^\x(x)\log\m_n^\x(x)\,d\n_{\check \m}d\n_\m
\le -\sum_{x \in B_{nL}^{\sf c}}\m_n(x)\log \m_n(x),
\]
by the Fubini theorem and the Jensen inequality.
By Lemma \ref{Lem:entropy_outside} \eqref{Eq:Lem:entropy1} and \eqref{Eq:Lem:entropy2},
for all large enough $L>4$ and for all positive integer $n$,
\begin{equation}\label{Eq:G3}
-\sum_{x \in B_{nL}^{\sf c}}\m_n(x)\log \m_n(x)
\le \e n D+C e^{-D' nL},
\end{equation}
where $C$, $D$ and $D'$ are positive constants independent of $n$ and $L$.

Combining \eqref{Eq:G1G2} and \eqref{Eq:G3},
we obtain for all $\e>0$ and for all $L, R$ with $\tilde \e_L<\e$, $\e_R<\e/L$ and for all large enough $n$,
\begin{align*}
\sup_{\m \in \Mcc} H(\m_n^\x) \le 3\e n D+\log \f_{n, R, L}+O(1),
\end{align*}
and thus together with \eqref{Eq:f}, we have $H(\m_n^\x)/n \to 0$ uniformly on $\m \in \Mcc$ as $n \to \infty$.
Since $H(\m_n^\x)=H(\m_n)-n h(\m)$ by \eqref{Eq:conditional_entropy},
and for each $n>0$, the $H(\m_n)$ is continuous in $\m \in \Mcc$ by Lemma \ref{Lem:entropy_outside},
we conclude that $\m \mapsto h(\m)$ is continuous on $\Mcc$.
\qed

\subsection*{Acknowledgment}
The author would like to thank J\'er\'emie Brieussel for bringing him the problem and discussions, and an anonymous referee for beneficial comments.
This work is partially supported by 
JSPS Grant-in-Aid for Scientific Research JP20K03602.

\bibliographystyle{alpha}
\bibliography{ns}

\begin{thebibliography}{GdlH90}

\bibitem[BB21]{BenjaminiBrieussel}
Itai Benjamini and J\'er\'emie Brieussel.
\newblock Noise sensitivity of random walks on groups.
\newblock arXiv:1901.03617v2, 2021.

\bibitem[BS00]{BonkSchramm}
M.~Bonk and O.~Schramm.
\newblock Embeddings of {G}romov hyperbolic spaces.
\newblock {\em Geom. Funct. Anal.}, 10(2):266--306, 2000.

\bibitem[Der80]{Derriennic80}
Y.~Derriennic.
\newblock Quelques applications du th\'eor\`eme ergodique sous-additif.
\newblock In {\em Conference on Random Walks (Kleebach, 1979)}, volume~74 of
  {\em Ast\'erisque}, pages 183--201, Paris, 1980. Soc. Math. France.

\bibitem[Der86]{Derriennic}
Y.~Derriennic.
\newblock Entropie, th\'{e}or\`emes limite et marches al\'{e}atoires.
\newblock In {\em Probability measures on groups, {VIII} ({O}berwolfach,
  1985)}, volume 1210 of {\em Lecture Notes in Math.}, pages 241--284.
  Springer, Berlin, 1986.

\bibitem[EK13]{ErschlerKaimanovich}
A.~Erschler and V.~A. Kaimanovich.
\newblock Continuity of asymptotic characteristics for random walks on
  hyperbolic groups.
\newblock {\em Funktsional. Anal. i Prilozhen.}, 47(2):84--89, 2013.

\bibitem[GdlH90]{GhysdelaHarpe}
\'{E}. Ghys and P.~de~la Harpe, editors.
\newblock {\em Sur les groupes hyperboliques d'apr\`es {M}ikhael {G}romov},
  volume~83 of {\em Progress in Mathematics}.
\newblock Birkh\"{a}user Boston, Inc., Boston, MA, 1990.
\newblock Papers from the Swiss Seminar on Hyperbolic Groups held in Bern,
  1988.

\bibitem[Gro87]{Gromovhyperbolic}
M.~Gromov.
\newblock Hyperbolic groups.
\newblock In S.~M. Gersten, editor, {\em Essays in group theory}, volume~8 of
  {\em Mathematical Sciences Research Institute Publications}, pages 75--263,
  New York, 1987. Springer-Verlag.

\bibitem[Kai00]{Kaimanovich-hyperbolic}
Vadim~A. Kaimanovich.
\newblock The {P}oisson formula for groups with hyperbolic properties.
\newblock {\em Ann. of Math. (2)}, 152(3):659--692, 2000.

\bibitem[Kal18]{KalaiICM2018}
Gil Kalai.
\newblock Three puzzles on mathematics, computation, and games.
\newblock In {\em Proceedings of the {I}nternational {C}ongress of
  {M}athematicians---{R}io de {J}aneiro 2018. {V}ol. {I}. {P}lenary lectures},
  pages 551--606. World Sci. Publ., Hackensack, NJ, 2018.

\bibitem[KV83]{KaimanovichVershik}
V.~A. Kaimanovich and A.~M. Vershik.
\newblock Random walks on discrete groups: boundary and entropy.
\newblock {\em Ann. Probab.}, 11(3):457--490, 1983.

\bibitem[Tan19]{Tdim}
Ryokichi Tanaka.
\newblock Dimension of harmonic measures in hyperbolic spaces.
\newblock {\em Ergodic Theory Dynam. Systems}, 39(2):474--499, 2019.

\end{thebibliography}

\end{document}